\documentclass{cimart}

\usepackage{faktor}
\usepackage{quiver}
\usepackage{tikz-cd}

\title{Description of DG-algebra structures via characters on a groupoid}

\author{Andronick Arutyunov and Oleg Muravev}

\authorinfo[A. Arutyunov]{Russian Academy of Sciences, Moscow, Russia}{andronick.arutyunov@gmail.com}

\authorinfo[O. Muravev]{Russian Academy of Sciences, Moscow, Russia}{oleg-muravev2001.mur@yandex.ru}

\abstract{%
   We construct a description of graded derivations in group algebras. Using this result for arbitrary grading of the group algebra, we describe all possible DG-algebra structures. Examples are given. The description is given in terms of characters on a groupoid analogous to the  groupoid of conjugate action.
   }



\VOLUME{33}
\YEAR{2025}
\ISSUE{3}
\NUMBER{13}
\DOI{https://doi.org/10.46298/cm.15775}

\begin{document}

\tableofcontents

\section{Introduction}

In several previous papers by A.S. Mishchenko and others \cite{arutyunov2017derivationsgroupalgebras, Arutyunov_2019, Alekseev_2023}, a method for studying the derivations of group algebras using the character space on the groupoid of conjugacy action was presented. In the present paper, a similar groupoid and characters on it are constructed in order to describe graded derivations on graded group algebras and all possible DG-algebra structures on group algebras. 

DG-algebras and graded derivations occur frequently. See e.g. \cite{Orlov_2023}, where they are studied in the context of algebraic geometry. Graded derivations in isolation are an even better-known subject, see e.g. \cite{Meirenken_2003}.
The classification of possible DG-structures has also been a focus of study \cite{Mao_2011, MAO2022389}.

The main object of our study will be the graded group algebra $\mathbb {C}[G]\cong \bigoplus A_{i}$ and its graded derivations, that is, linear maps $d: \mathbb {C}[G] \rightarrow \mathbb {C}[G]$ satisfying the graded Leibniz rule 
$$d(uv) = d(u)v + (-1)^{|u|}ud(v), \quad \forall u,v \in \cup A_{i}.$$ 

In Section \ref{sec-grouppoid}, we present the structure of a groupoid which is connected to the graded algebra structure. We will give a description of graded derivations as locally finite characters on this groupoid (see Theorem \ref{Th-d(g)=}).
In Section \ref{examples}, we describe the inner, quasi-inner, and central derivations in the graded case. Using this construction in the case of DG-algebras, in Section \ref{DG-struc}, we obtain a description of DG-algebra structures in terms of characters and an isomorphism criterion for DG-algebras. In addition, at the end of Section \ref{DG-struc}, we give several examples of gradings on group algebras and derivations that define a DG-algebra structure on them. The definitions and properties from combinatorial group theory that are necessary to understand the present paper can be found in \cite{Lyndon_2001}.  
\section{The groupoid, its characters, and graded derivations}
\label{sec-grouppoid}

Let's start with the definitions and fix the notation for future reference. We will study derivations of group algebras with a DG-structure. The standard definition of a group algebra is

\begin{definition}
Let $G$ be a group and let $K$ be a field. The group algebra $K[G]$ is the associative algebra over $K$ the elements of which are all possible formal finite sums of the form $\sum_{g \in G}k_{g}g,$ $g \in G,$ $k_{g} \in K$, and the operations are defined by: $$\sum_{g \in G}a_{g}g + \sum_{g \in G}b_{g}g = \sum_{g \in G}(a_{g} + b_{g})g,$$
$$(\sum_{g \in G}a_{g}g)(\sum_{g \in G}b_{g}g) = \sum_{h \in G}\Big(\sum_{xy = h \in G}(a_{x}b_{y})h\Big).$$
The elements of the group $G$ form a basis of $K[G]$. Multiplication of the basis elements in the group algebra is induced by the group multiplication.
\end{definition}

Furthermore, we need the notion of a groupoid. Usually, it is defined as a category in which every morphism is invertible. We will work with the following construction associated with a finitely generated group $G$, so this category is small (the classes of objects and morphisms are sets).

\begin{definition}\label{def-gruppoid}
A groupoid $\Gamma$ has a set of objects $\operatorname{Obj}(\Gamma)$ and a set of morphisms $\operatorname{Hom}(\Gamma)$. Each morphism has a source $s(\phi)$ and a target $t(\phi)$, i.e. each morphism is in the set $\operatorname{Hom}(s(\phi), t(\phi))$. Composable morphisms have an associative operation (composition). The following conditions hold:
\begin{enumerate}
    \item For each object $a\in \operatorname{Obj}(\Gamma)$, there is a neutral endomorphism $1_a\in \operatorname{Hom}(a,a)$ such that $\forall \phi\in \operatorname{Hom}(b,a)$: $1_a\circ \phi = \phi$ and $\forall \phi\in \operatorname{Hom}(a,b)$ $ \phi\circ 1_a = \phi$.
    \item Each morphism has an inverse: for every $\phi\in \operatorname{Hom}(a,b)$ there exists $\psi\in \operatorname{Hom}(b,a)$ such that $\phi\circ \psi = 1_b$ and $\psi\circ\phi = 1_a$.
\end{enumerate}
\end{definition}

\subsection{Groupoid of conjugacy action}
The following construction generalizes the one from \cite{arutyunov2017derivationsgroupalgebras}.
Suppose we have a $\mathbb {C}$--graded group algebra $\mathbb {C}[G]\cong \bigoplus A_{i}$. The group $G$ is assumed to be finitely generated.
We will define the groupoid of conjugacy action $\Gamma$ for the group $G$ (in the sense of Definition \ref{def-gruppoid}) as follows:
$$
\Gamma =\begin{cases}
\operatorname{Obj}(\Gamma) = \{\pm g \mid g \in G\},\\
\operatorname{Hom}(\Gamma) = \{(u,v) \in \pm G \times  \pm G\ \mid s((u,v)) = (-1)^{|v|}v^{-1}u,\;  t((u,v)) = uv^{-1}\}.\\
\end{cases}
$$
Consider two morphisms $\phi = (u_{1},v_{1})$ and $\psi = (u_{2},v_{2})$ such that $t(\phi) = s(\psi)$, that is, the composition $\psi \circ \phi$ should be defined. Define it as follows: 
$$(u_{2} , v_{2})\circ(u_{1} , v_{1}) := \big((-1)^{|v_{2}|}v_{2} u_{1}, v_{2}v_{1}\big).$$

This formula is conveniently depicted as a diagram.
 
 \[\begin{tikzcd}
  & {(-1)^{|v_1|}v_1^{-1}u_1 = u_2v_2^{-1}} \\
  {u_1v_1^{-1}} && {(-1)^{|v_2|}v_2^{-1}u_2}
  \arrow["{(u_2,v_2)}"{pos=0.7}, from=1-2, to=2-3]
  \arrow["{(u_1,v_1)}"{pos=0.3}, from=2-1, to=1-2]
  \arrow["{(u_2,v_2)\circ (u_1,v_1)}"', from=2-1, to=2-3]
\end{tikzcd}\]
 
\begin{proposition} \label{t1} % сразу указываем ссылку
$\Gamma$ is a groupoid. \end{proposition}

\begin{proof}
First of all, we note that $\exists g \rightarrow h \Leftrightarrow \exists t: h = (-1)^{|t|}tgt^{-1}.$

Consider $\phi = (u_1,v_1)$ and $\psi = (u_2,v_2)$ such that $t(\phi) = s(\psi)$. Define $s(\phi) = a$, $t(\phi) = s(\psi) = b$ and $t(\psi) = c$.

We need to check that
\begin{center}
$s(\psi\circ\phi) = a$ and $t(\psi\circ\phi) = c$.
\end{center}

By definition,
\begin{center}
$s(\psi\circ\phi) = (-1)^{|v_{2}v_{1}|\mathstrut}v_{1}^{-1\mathstrut}v_{2}^{-1\mathstrut}(-1)^{|v_{2}|\mathstrut}v_{2}u_{1}$ = $(-1)^{|v_{1}|\mathstrut}v_{1}^{-1\mathstrut}u_{1} = a$.
\end{center}
Calculate that 
\begin{center} 
$t(\psi\circ\phi)$ = $(-1)^{|v_{2}|\mathstrut}v_{2}u_{1}v_{1}^{-1\mathstrut}v_{2}^{-1\mathstrut} = u_{2}v_{2}^{-1\mathstrut} = c$, because $u_{1}v_{1}^{-1\mathstrut} = (-1)^{|v_{2}|\mathstrut}v_{2}^{-1\mathstrut}u_{2}.$
\end{center}
Here we shall prove that the composition is associative. Let $\phi_{1},\phi_{2},\phi_{3}$ be morphisms such that $s(\phi_{1}) = t(\phi_{2})$ and $s(\phi_{2}) = t(\phi_{3})$. We need to show that
\begin{center}
$\phi_{1}\circ(\phi_{2}\circ\phi_{3}) = (\phi_{1}\circ\phi_{2})\circ\phi_{3}.$
\end{center}
We have proved above that these compositions are well-defined. The left side is
\begin{center}
$\phi_{2}\circ\phi_{3} = ((-1)^{|v_{2}|}v_{2} u_{3}, v_{2}v_{3}))$
$\phi_{1}\circ(\phi_{2}\circ\phi_{3}) = (u_{1},v_{1})\circ((-1)^{|v_{2}|}v_{2} u_{3}, v_{2}v_{3})) = ((-1)^{|v_{1}|+|v_{2}|}v_{1}v_{2}u_{3},v_{1}v_{2}v_{3}).$
\end{center}
While the right one is equal to it,
\begin{center}
$(\phi_{1}\circ\phi_{2})\circ\phi_{3} = ((-1)^{|v_{1}|}v_{1}u_{2},v_{1}v_{2})\circ(u_{3},v_{3}) = ((-1)^{|v_{1}|+|v_{2}|}v_{1}v_{2}u_{3},v_{1}v_{2}v_{3}),$
\end{center}
proving associativity.
\\Let us show that there is an element $1_{g}$ such that
\begin{center}
$s(1_{g}) = t(1_{g}) = g$,\\
$\phi\cdot1_{s(\phi)} = \phi = 1_{t(\phi)}\cdot\phi$, $\forall\phi\in$ $\operatorname{Hom}(\Gamma).$
\end{center}

%Диаграмма для обратного и единичного
Take $1_{g} = (g,e)$, where $e$ is the identity of the group $G:$
\[\begin{tikzcd}
	(-1)^{|e|}e^{-1}g = g = ge^{-1}
	\arrow["{(g,e)}", from=1-1, to=1-1, loop, in=55, out=125, distance=10mm]
\end{tikzcd}\]

Indeed, choosing an arbitrary morphism $(u,v)$, we have:
\begin{center}
$(uv^{-1},e)(u,v) = (u,v)$;\\[1mm]
$(u,v)((-1)^{|v|}v^{-1}u,e) = ((-1)^{|v|}v(-1)^{|v|}v^{-1}u,e) = (u,v)$.
\end{center}
We have now shown that $\Gamma$ is a category. It remains to check the invertibility of the arrows.
\\We claim the inverse of a morphism of the form $(u,v)$ is $((-1)^{|v|}v^{-1}uv^{-1}, v^{-1})$.

Indeed, 
\begin{center}
$(u,v)((-1)^{|v|}v^{-1}uv^{-1}, v^{-1}) = (uv^{-1}, e) = (t(\phi),e),$\\[1mm]
$((-1)^{|v|}v^{-1}uv^{-1}, v^{-1})(u,v) = ((-1)^{|v|}v^{-1}u,e) = (s(\phi),e).$
\end{center}

\[\begin{tikzcd}
	{(-1)^{|v|}v^{-1}u} && {uv^{-1}}
	\arrow["{(u,v)}", shift left=2, curve={height=-18pt}, from=1-1, to=1-3]
	\arrow["{((-1)^{|v|}v^{-1}uv^{-1}, v^{-1})}", shift left=3, curve={height=-18pt}, from=1-3, to=1-1]
\end{tikzcd}\]

This proves Proposition  \ref{t1}. 
\end{proof} 

The following properties of the groupoid are easily verified by direct calculation.

\begin{proposition}

\begin{enumerate} 
\item The endomorphisms $\operatorname{Hom}(a,a), \ a \in\pm G$, admit the following description:
\begin{center}
$\operatorname{Hom}(a,a) = \{\chi \in \operatorname{Hom}(\Gamma) \mid s(\chi) = (-1)^{|t|}t^{-1}a,$ $t(\chi) = at^{-1},$ $t \in Z(a)$ and $t \in$ $A_{2i} \}.$    
\end{center}

\item For $[a] := \{b\in \pm G \mid \exists t\in \pm G: b = (-1)^{|t|}t^{-1}at\}$, denote by $\Gamma_{[a]}$ the subgroupoid such that $\operatorname{Obj}(\Gamma_{[a]}) = [a]$. The following properties hold:
\begin{enumerate}
    \item $\operatorname{Hom}(\Gamma_{[a]}) = \{ \phi \in \operatorname{Hom}(\Gamma) \mid s(\phi) \in [a], t(\phi) \in [a] \}$.
    \item $\Gamma = \coprod \Gamma_{[a]}$.
\end{enumerate}

\item The left action $\operatorname{Hom}(a,a)\times \operatorname{Hom}(a,b) \rightarrow \operatorname{Hom}(a,b)$: $(\phi,\psi) \to \psi \circ \phi$ is free and transitive. 
\end{enumerate} 
\end{proposition}

Let us show that, with the help of characters of $\Gamma$, it is possible to describe derivations satisfying the graded Leibniz rule.
\begin{definition}

The linear operator $d$: $\mathbb {C}[G]$ $\to$ $\mathbb {C}[G]$ satisfies the graded Leibniz rule if 
\begin{center}
$d(uv) = d(u)v + (-1)^{|u|}ud(v)$, $\forall \  u,v \in \cup A_{i}.$
\end{center}
\end{definition}

\begin{definition}
The map $\chi$ : $\operatorname{Hom}(\Gamma)$ $\to$ $\mathbb {C}$ is a character if it satisfies the identity
 \[
 \chi(\psi\circ\phi)=\chi(\psi)+\chi(\phi):
 \qquad
 \begin{tikzcd}
  && \bullet \\
  \bullet &&& {} & \bullet
  \arrow["{\chi(\psi)}", from=1-3, to=2-5]
  \arrow["{\chi(\phi)}", from=2-1, to=1-3]
  \arrow["{\chi(\psi\circ\phi)=\chi(\psi)+\chi(\phi)}"', from=2-1, to=2-5]
\end{tikzcd}\]
\end{definition}  

\begin{definition}
    A character $\chi$ is locally finite if for all $v \in \pm G$ we have $\chi(u,v) = 0$ for all except a finite number of $u$.
\end{definition}

We will denote by $\chi(\Gamma)$  the space of locally finite characters. In  Section \ref{characters} we will show that $\chi(\Gamma)$ has a natural Lie algebra structure.

The graded Leibniz rule in a group algebra can be expressed in terms of  characters on the groupoid $\Gamma$ as follows.
 
\begin{proposition} \label{t2} % сразу указываем ссылку
For any graded derivation $d$, there exists a locally finite character $\chi$ such that
\begin{equation}
\label{1}
d(u) = \sum\limits_{g\in G}\Big(\sum\limits_{h\in \pm G}\chi(\phi(h,g)\lambda^{h}\Big)g, \  \forall u \in A_{i}.
\end{equation}
\end{proposition} 

\begin{proof}
It is enough to check the statement on the elements $g\in \pm G$, since they generate $\mathbb {C}[G]$ as a vector space. 

For any element $g$ of group G we have $d(g) = \sum\limits_{h\in G}d_{g}^{h}h$.

Let us define $\chi(\phi(h,g)):=d_{g}^{h}$, $\chi(\phi(-h,g)):= -d_{g}^{h}$ and show that $\chi$ it is a character.
We have 
\begin{center}
$d(g_{2}g_{1}) = d(g_{2})g_{1} + (-1)^{|g_{2}|}g_{2}d(g_{1}),$
\end{center}
or after a basis expansion
\begin{center}
$\sum\limits_{h\in G}d_{g_{2}g_{1}}^{h}h = \sum\limits_{h\in G}d_{g_{2}}^{h}hg_{1} + (-1)^{|g_{2}|}g_{2}\sum\limits_{h\in G}d_{g_{1}}^{h}h.$
\end{center}

It is clear that summing over all $h$ and all $gh$ is equivalent, so we can convert the right-hand side to the following form:
\begin{equation}
\label{2}
\sum\limits_{h\in G}d_{g_{2}g_{1}}^{h}h = \sum\limits_{h\in G}d_{g_{2}}^{hg_{1}^{-1}}h+ \sum\limits_{h\in  G}d_{g_{1}}^{(-1)^{|g_{2}|}g_{2}^{-1}h}h \Rightarrow d_{g_{2}g_{1}}^{h} = d_{g_{2}}^{hg_{1}^{-1}} + d_{g_{1}}^{(-1)^{|g_{2}|}g_{2}^{-1}h}.
\end{equation}

Let's make a substitution:
\begin{center}
$h_{1} = hg_{1}^{-1}, \ \ h_{2} = (-1)^{|g_{2}|}g_{2}^{-1}h.$	
\end{center}

Then the last equality (\ref{2}) takes the form:
\begin{center}
$d_{g_{2}g_{1}}^{(-1)^{|g_{2}|}g_{2}h_{2}} = d_{g_{2}}^{h_{1}} + d_{g_{1}}^{h_{2}}$ $\Leftrightarrow$ $\chi(\phi((-1)^{|g_{2}|}g_{2}h_{2},g_{2}g_{1})) = \chi(\phi(h_{1},g_{2})) + \chi(\phi(h_{2},g_{1})).$
\end{center}

It remains only to note that $((-1)^{|g_{2}|}g_{2}h_{2},g_{2}g_{1}) = (h_{1},g_{2})(h_{2},g_{1}).$ 
\end{proof} A corollary of the above statement is the following theorem, proved similarly to the non-graded case, see Theorem 1 of \cite{Arutyunov_2020}.
\begin{theorem}\label{Th-d(g)=}
For any graded derivation $d$ there exists a unique locally finite character $\chi_{d}$ satisfying

\begin{center}
$d(a) = (-1)^{|a|}a\Big(\sum\limits_{t\in G}\chi((-1)^{|a|}at, a)t\Big).$
\end{center}
\end{theorem}
There is therefore an isomorphism $\operatorname{Der}(\mathbb {C}[G]) \cong \chi(\Gamma)$. 

\begin{proof}
To prove this, we need to apply formula (1) to $a \in \pm G$, write an arbitrary element $g\in G$ as $g = at$, and then drop the terms with a zero coefficient from the sum. See the proof of Theorem 1 in \cite{Arutyunov_2020} for more details. 
\end{proof}

\subsection{Characters}
\label{characters}
Locally finite characters are isomorphic to derivations, which means that the commutator of derivations yields a binary operation on characters. Formula (\ref{1}) allows us to define a commutator in character space. Let $d_{1}, d_{2}$ be the derivations defined by the characters $\chi_{d_{1}}, \chi_{d_{2}}$. Denote the character that corresponding to the commutator $[d_{1}, d_{2}]$ by $\chi_{[d_{1}, d_{2}]}$. This defines an operation on the space of characters. 
\begin{equation}
\label{3}
\{\chi_{d_{1}}, \chi_{d_{2}}\} := \chi_{[d_{1}, d_{2}]}.
\end{equation}
For this operation, we have the following equality:
\begin{proposition}
The operation on characters above can be expressed as
\begin{equation}
\label{4}
\{\chi_{d_{1}}, \chi_{d_{2}}\}(a,g) = \sum\limits_{h \in G}\Big(\chi_{d_{1}}(a,h)\chi_{d_{2}}(h,g) - \chi_{d_{2}}(a,h)\chi_{d_{1}}(h,g)\Big).
\end{equation}
\end{proposition}
\begin{proof}
Let $g \in G$. Then we have 
\begin{align*}
d_{1}(g) &= \sum\limits_{h \in  G}\chi_{d_{1}}(h,g)h,\\
d_{2}(g) &= \sum\limits_{h \in  G}\chi_{d_{2}}(h,g)h,\\
\{\chi_{d_{1}}, \chi_{d_{2}}\}(g) &= \sum\limits_{h \in G}\{\chi_{d_{1}}, \chi_{d_{2}}\}(a,g)a.
\end{align*}	
Write down the equation for the commutator
\begin{center}
$[d_{1},d_{2}] = d_{1}d_{2} - d_{2}d_{1}$,\\
$d_{1}d_{2}(g) = \sum\limits_{h \in G}\chi_{d_{2}}(h,g)\Big(\sum\limits_{a \in \pm G}\chi_{d_{1}}(a,h)a\Big)$,\\
$d_{2}d_{1}(g) = \sum\limits_{h \in G}\chi_{d_{1}}(h,g)\Big(\sum\limits_{a \in \pm G}\chi_{d_{2}}(a,h)a\Big)$.
\end{center}
Thanks to local finiteness, we can swap the sums in the last equation and get 
\begin{center}
$[d_{1},d_{2}](h) = \sum\limits_{a \in G}\Big(\sum\limits_{h \in G}\chi_{d_{2}}(h,g)\chi_{d_{1}}(a,h) - \chi_{d_{1}}(h,g)\chi_{d_{2}}(a,h)\Big)a$.
\end{center}
The value of $\{\chi_{d_{1}}, \chi_{d_{2}}\}(a,g)$ is the coefficient of $a$, so
\begin{center}
$\{\chi_{d_{1}}, \chi_{d_{2}}\}(a,g) = \sum\limits_{h \in G}\chi_{d_{1}}(a,h)\chi_{d_{2}}(h,g) - \chi_{d_{2}}(a,h)\chi_{d_{1}}(h,g)$.
\end{center} \end{proof}

It follows that $(\chi_{\Gamma}, \{\cdot,\cdot \})$ is a Lie algebra, isomorphic to the algebra of all graded derivations by Proposition \ref{t2}.

\section[Examples of graded derivations and their description]{Examples of graded derivations and their description in terms of characters}
\label{examples}
We will illustrate the correspondence between characters on a groupoid and graded derivations using examples of inner derivations, central derivations, and quasi-inner derivations. They will serve as easy examples of such operators.

\subsection{Inner derivations}

Consider an element $a \in \pm G$. Define the mapping $\chi^{a}: \operatorname{Hom}(\Gamma) \rightarrow \mathbb{C}$ as follows (assuming $a \ne b$):
\begin{equation}
\label{eq-charinner}
\chi^{a}(\phi) := 
 \begin{cases}
   1, &\text{$\phi \in \operatorname{Hom}(b,a)$}\\
   -1, &\text{$\phi \in \operatorname{Hom}(a,b)$}\\
   0, &\text{otherwise.}
 \end{cases}
\end{equation}
Clearly $\chi^{a}$  is trivial on endomorphisms.
\begin{proposition}\label{Inn}
The map $\chi^{a}$ defines an inner derivations for $a \in \pm G$ by the formula 
\begin{equation}
\label{InnForm}
d_{a} : x \mapsto [a,x]_{grad} = ax - (-1)^{|x|}xa.
\end{equation}
\end{proposition}
\begin{proof}
Let's check that $d_{a}$ is a graded derivation (that is, satisfies the graded Leibniz rule):
\begin{center}
$d_{a}(x)y + (-1)^{|x|}xd_{a}(y) = (ax - (-1)^{|x|}xa)y + (-1)^{|x|}x(ay - (-1)^{|y|}ya) = axy - (-1)^{|xy|}xya = d_{a}(xy)$.
\end{center}
Let's check that it really satisfies the formula \eqref{InnForm}. Due to linearity, it is enough to check the statement for the element $g \in \pm G$. According to the previous corollary, we have an equality:
\begin{center}
$d_{a}(g) = \sum\limits_{h\in G}\chi(h, g)h.$
\end{center}
We have that $\chi(h, g) \ne 0 \Leftrightarrow t((h, g)) = a$ or $s((h, g)) = a$. Let's consider these two cases:
\begin{enumerate}
\item $s((h, g)) = a \Leftrightarrow (-1)^{|g|}g^{-1}h = a \Leftrightarrow h = (-1)^{|g|}ga.$
\item $t((h, g)) = a \Leftrightarrow hg^{-1} = a \Leftrightarrow h = ag.$
\end{enumerate}
Therefore, we have: $d_{a}(g) = ag - (-1)^{|g|}ga.$
\end{proof}

\subsection{Central derivations}
We will now give an example of a derivation that is not inner. We will show that it admits a simple explicit form and is not trivial on endomorphisms. That entails that it cannot be expressed as a (perhaps formal) sum of inner derivations.
\begin{definition}
That the map $\tau: G \rightarrow \mathbb{C}$ is a graded group character if it satisfies the following property:
\begin{center}
$\tau(ab) = \tau(a) + (-1)^{|a|}\tau(b),$  $\forall a, b \in G$.
\end{center}
\end{definition}
Let's fix a graded group character $\tau$, and let $z \in Z(G)$ be a central element of the group.  
\begin{definition}
Define the map $d_{z}^{\tau}: \mathbb{C}[G] \rightarrow \mathbb{C}[G]$ on the generators $g \in \pm G$ as
\begin{center}
$d_{z}^{\tau}: g \rightarrow \tau(g)gz$
\end{center}
and extend it by linearity to the whole algebra $\mathbb{C}[G]$. We will call it a central graded derivation.
\end{definition}
\begin{proposition}
The map $d_{z}^{\tau}$ is a graded derivation.
\end{proposition}
\begin{proof}
It is enough to check the graded Leibniz rule on arbitrary generators $u, v \in G$:
\begin{center}
$d_{z}^{\tau}(uv) = \tau(uv)uvz = \tau(u)uzv + (-1)^{|u|}u\tau(v)vz = d_{z}^{\tau}(u)v + (-1)^{|u|}ud_{z}^{\tau}(v)$.
\end{center}
On arbitrary elements of a group algebra, the graded Leibniz rule is satisfied due to the linearity of $d_{z}^{\tau}$.
\end{proof}
\begin{proposition}\label{prop-centralnoninner}
Nontrivial central derivations are not inner.
\end{proposition}
\begin{proof}
Let $\chi_{\tau, z}$ be a character that corresponds to the derivation $d_{z}^{\tau}$. By Theorem \ref{Th-d(g)=} we have the equality
\begin{center}
$d_{z}^{\tau}(g) = \tau(g)gz = g\Big(\sum\limits_{t\in G}\chi((-1)^{|g|}gt, g)(-1)^{|g|}t\Big)$.
\end{center}
Therefore, this character is not equal to 0 only in the case $t = z$:
\begin{center}
$\chi((-1)^{|g|}gz, g) = (-1)^{|g|}\tau(g)$.
\end{center}
Define the morphism $\phi = ((-1)^{|g|}gz, g)$ and calculate its beginning and
\begin{center}
$s(\phi) = z,\ \
t(\phi) = (-1)^{|g|}z.$
\end{center}
Then in the case when $g \in A_{2n}$, we have: $s(\phi) = t(\phi)$ $\Rightarrow \chi_{\tau, z}$ is non--trivial on a loop $\phi$, thus $d_{z}^{\tau}$ it is not an inner derivation.
\end{proof}

A direct check shows that inner and central derivations constitute subalgebras. For more details, see \cite{Arutyunov_2020}.

\subsection{Quasi-inner derivations}

As seen from the formula \eqref{eq-charinner}, inner derivations are trivial on endomorphisms; however, the converse is not true. We will call such derivations quasi-inner.

\begin{definition}
Define $\operatorname{QInnDer}(\mathbb{C}[G])$ as follows:
\begin{equation}
\operatorname{QInnDer}(\mathbb{C}[G]) = \left\{ d \in \operatorname{Der}     \mid \forall a\in \operatorname{Obj}(\Gamma), \forall \phi\in \operatorname{Hom}(a,a) : \chi_{d}(\phi) = 0 \right\}.
\end{equation}
\end{definition}

Quasi-inner derivations are clearly a subspace, but the more significant statement is that they constitute an ideal in the algebra of all derivations.

\begin{proposition}
$\operatorname{QInnDer}(\mathbb{C}[G])$ $\triangleleft$ $\operatorname{Der}(\mathbb{C}[G])$ is an ideal.
\end{proposition}

\begin{proof}
Let us prove that $\operatorname{QInnDer}(\mathbb{C}[G])$ is a subalgebra, i.e.
\begin{center} 
$d_{1}$, $d_{2}$ $\in$ $\operatorname{QInnDer}(\mathbb{C}[G]$) $\Rightarrow$ $\left[ d_{1}, d_{2}\right] \in \operatorname{QInnDer}(\mathbb{C}[G])$.
\end{center}
The character $\chi_{d_{1}}$ can be represented as
\begin{center}
$\chi_{d_{1}} = \sum\limits_{a\in G}\lambda^{a}\chi^{a}$, $\lambda^{a} \in \mathbb{C}$,
\end{center}
where $\chi_{d_{1}}$ defined by formula (\ref{eq-charinner}), since $\chi^{a}$ is trivial on endomorphisms.
Similarly for the character $\chi_{d_{2}}$:
\begin{center}
$\chi_{d_{2}} = \sum\limits_{b\in G}\lambda^{b}\chi^{b}$, $\lambda^{b} \in \mathbb{C}$.
\end{center}
Because of the bilinearity of the commutator,
\begin{center}
	$\{\chi_{d_1}, \chi_{d_2}\} = \sum\limits_{a \in G}\sum\limits_{b \in G}\lambda^a\mu^b \{\chi^a, \chi^b\}$.
\end{center}

The commutator of the $d^a$ and $d^b$ can be represented as follows:
\begin{center}
	$[d^a, d^b] = d^{ab} - d^{ba}$.
\end{center}

Then define the character $\{\chi^a, \chi^b\}$ by the following formula
\begin{center}
    $\{\chi^a, \chi^b\} = \chi^{ab} - \chi^{ba}$.
\end{center}

Claim:
\begin{center}
    $\{\chi_{d_1}, \chi_{d_2}\} = \sum\limits_{a \in G}\sum\limits_{b \in G}\lambda^a\mu^b \chi^{ab} - \sum\limits_{a \in G}\sum\limits_{b \in G}\lambda^a\mu^b \chi^{ba}$,
\end{center}

thus $\{\chi_{d_1}, \chi_{d_2}\} \in \operatorname{QInnDer}(\mathbb{C}[G])$. To prove that, consider the value of the character on the loop $(uz, z), z \in Z_G(u)$:
\begin{center}
	$\{\chi_{d_1}, \chi_{d_2}\}(uz, z) = \sum\limits_{ab = zuz^{-1}}\lambda^a\mu^b - \sum\limits_{ab = u}\lambda^a\mu^b - \sum\limits_{ba = u}\lambda^a\mu^b + \sum\limits_{ba = zuz^{-1}}\lambda^a\mu^b = 0$.
\end{center}

We will now turn to the proof of the theorem. Let's represent $\chi_{d_0}$ as
\begin{center}
    $\chi_{d_0} = \sum\limits_{a \in G} \lambda^a\chi^a$.
\end{center}
 
	It is sufficient to prove the theorem for $\chi^a$ and extend the result to the $\chi_{d_0}$ using bilinearity of the commutator. Consider the character $\{\chi_d, \chi^a\}$. We need to prove that it is trivial on endomorphisms, i.e., that for all $b \in G$ and for $z \in Z_G(b)$ it satisfies $\{\chi_d, \chi^a\}(bz, z) = 0$. By Proposition \ref{t1}:
    \begin{center}
    	$\{\chi_d, \chi^a\}(bz, z) = \sum\limits_{h\in G} \chi_d(bz, h)\chi^a(h, z) - \chi^a(bz, h)\chi_d(h, z)$.
    \end{center}

Here $\chi^a(h, z) \neq 0$ only in two cases: when $h = za$ and when $h = az$; and $\chi^a(bz, h) \neq 0$ only in two cases: when $h = bza^{-1}$ and when $h = a^{-1}bz$. This means that

	\begin{center}
    	$\{\chi_d, \chi^a\}(bz, z) = \chi_d(bz, za) - \chi_d(bz, az) + \chi_d(a^{-1}bz, z) - \chi_d(bza^{-1}, z)$.
    \end{center}

But also
\begin{center}
	$(a^{-1}bz, z) \circ (bz, za) = (bz, az) \circ (bza^{-1}, z)$.
\end{center}
Hence 
\begin{center}
	$\chi_d(bz, za) + \chi_d(a^{-1}bz, z) = \chi_d(bz, az) + \chi_d(bza^{-1}, z)$,
    \\$\{\chi_d, \chi^a\}(bz, z) = 0$.
\end{center}

\end{proof}

It is easy to see that central derivations are not quasi-inner by repeating the reasoning from Proposition \ref{prop-centralnoninner}.

\section[Description of the structure of a DG-algebra in terms of characters]{Description of the structure of a DG-algebra in terms of characters on a groupoid}
\label{DG-struc}

Our next goal is to describe the possible DG-algebra structures up to isomorphism.
We follow the definitions of \cite{F_lix_2001}.

\begin{definition}
Let $A\cong \bigoplus A_{i}$ be a $\mathbb{Z}$--graded algebra equipped with a derivation
\begin{center}
$d: A \rightarrow A$
\end{center}
of degree $1$ (the case of a cochain complex) or $-1$ (the case of a chain complex) and satisfying the  following conditions:
\begin{enumerate} 
\item $d \circ d = 0$,
\item $d(uv) = d(u)v + (-1)^{|u|}ud(v)$.
\end{enumerate}
We say that $(A, d)$ is a differential graded algebra or DG-algebra.       
\end{definition}

We will give the necessary and sufficient conditions from the point of view of groupoid characters for a derivation to define a DG-algebra structure on a graded group algebra. We will describe the cochain complex case, since the reasoning in the chain complex is similar.
\begin{theorem}
\label{DG_iff}
A character $\chi$ on the groupoid $\Gamma$ defines a DG-algebra structure if and only if the following conditions are satisfied:

\begin{enumerate}
\item The convolution is trivial: 
\begin{equation}
\label{1_Th-3.1}
\sum\limits_{h \in G}\chi(h, g)\chi(x,h) = 0, \quad \forall g, x \in G.
\end{equation} 
\item For all $(h,g) \in Hom(\Gamma_{A_{i}})$, $i \neq 1$, we have $\chi(h, g) = 0$. 

\end{enumerate}
\end{theorem}
The main idea is to prove that the first condition is equivalent to $d^2(g) = 0$ and the second one is equivalent to $d(A^i) \subseteq A^{i+1}$.

\begin{proof} 
Using Proposition \ref{t2}, we rewrite the defining properties of a DG-algebra in terms of characters.
\begin{enumerate}
\item By direct calculation we can see that $$d^{2}(g) = d\Big(\sum\limits_{h \in G}d^{h}_{g}h\Big) = \sum\limits_{h \in G}d^{h}_{g}\Big(\sum\limits_{x \in G}d^{x}_{h}x\Big) =$$$$ \sum\limits_{h \in G}(h, g)\Big(\sum\limits_{x \in G}(x, h)x\Big) = \sum\limits_{x \in G}\Big(\sum\limits_{h \in G}(\chi(h,g)\chi(x,h))\Big)x.$$

\item Let $g\in A^{i}$, then $d(g) = \sum\limits_{h \in G}\chi(h, g)h \in A^{i+1} \Leftrightarrow \chi(h, g) = 0$ when $h\notin A^{i+1}$.\\
Note that if $u \in A^{i+1}, v \in A^{i}$, then for $\phi = (u,v)$ we have: $s(\phi) = (-1)^{|v|}v^{-1}u \in A^{1}$, $t(\phi) = uv^{-1} \in A^{1}$, which means that in the case of a DG-algebra, the characters may not be equal to $0$ only on morphisms between elements from $A^{1}$.
\end{enumerate}
\end{proof}

Now, let us recall the definition of an isomorphism of DG-algebras.
\begin{definition}
\label{DG_def}
An isomorphism $f$: $A_{1}$ $\to$ $A_{2}$ between two DG-algebras $(A_{1}, d_{1})$ and $(A_{2}, d_{2})$ is an algebra isomorphism such that following properties hold:
\begin{enumerate}
\item f respects the grading, $f(A^{i}_{1}) \subset A^{i}_{2}$;
\item for all $a \in A_{1}$, $f(d_{1}(a)) = d_{2}(f(a))$.
\end{enumerate}
\end{definition}
Now we describe the structures of DG-algebras up to isomorphism per Definition \ref{DG_def}.
%Suppose now there is a situation that we are given a group algebra $\mathbb{C}[G]$, we want to classify the DG--structures on it up to isomorphism, in terms of characters on the groupoid. The following theorem holds

Let $X_{i}, i = 1,2$ be two matrices (infinity ones, if $G$ is an infinite group) in following sense: 
$$X_{i}: G \times G \to \mathbb{C},$$ $$(x,h) \mapsto \chi_{i}(x,h).$$
\begin{theorem}
Two DG-algebras $(\mathbb{C}[G], d_{1})$ and $(\mathbb{C}[G], d_{2})$ are isomorphic iff there is a matrix $C$ such that 
\begin{enumerate}
\item matrices $X_{1}$ and $X_{2}$ are conjugate, $CX_{2} = X_{1}C;$
\item the mapping f is in $\operatorname{Aut}(\mathbb{C}[G])$.
\end{enumerate}

\end{theorem}

In other words, that means that sets of constants $c_{h}^{g}$, gives $f(g) = \sum\limits_{h \in G}c_{h}^{g}h \in \operatorname{Aut}(\mathbb{C}[G])$ and
\begin{equation}
\label{equality}
\sum\limits_{h \in G}c_{h}^{g}\chi_{2}(x, h) = \sum\limits_{h \in G}\chi_{1}(h, g)c_{x}^{h}, \forall g, x \in G,
\end{equation}
where $\chi_{1}$ and $\chi_{2}$ denotes the characters that define the derivations $d_{1}$ and $d_{2}$.

\begin{proof}
Suppose that we have an isomorphism of DG-algebras $f: (\mathbb{C}[G], d_{1}) \rightarrow (\mathbb{C}[G], d_{2})$, given on generators by the formula $f(g) = \sum\limits_{h \in G}c_{h}^{g}h$. By definition, we have
\begin{center}
$f(d_{1}(g)) = d_{2}(f(g))$.
\end{center}
Then
\begin{center}
$f(d_{1}(g)) = f\Big(\sum\limits_{h \in G}\chi_{1}(h,g)h\Big) = \sum\limits_{h \in G}\chi_{1}(h,g)f(h) = \sum\limits_{h \in G}\chi_{1}(h,g)\big(\sum\limits_{x \in G}c_{x}^{h}x\Big) = \Big(\sum\limits_{h \in G}\chi_{1}(h,g)c_{x}^{h}\Big)x$;
\end{center}
\begin{center}
$d_{2}(f(g)) = d_{2}\Big(\sum\limits_{h \in G}c_{h}^{g}h\Big) = \sum\limits_{h \in G}c_{h}^{g}d_{2}(h) = \sum\limits_{h \in G}c_{h}^{g}(\sum\limits_{x \in G}\chi_{2}(x,h)x) = \Big(\sum\limits_{h \in G}c_{h}^{g}\chi_{2}(x,h)\Big)x$.
\end{center}
Whence we get our equality \eqref{equality}.
To prove the converse, just define $f(g) = \sum\limits_{h \in G}c_{h}^{g}h$ and for the same reasons, get the necessary equalities.
\end{proof}

\begin{corollary}
In the case when $f \in \operatorname{Aut(G)}$ the equality \eqref{equality} takes the form
\begin{center}
$\chi_{1}(h,g) = \chi_{2}(f(h),f(g)).$
\end{center}
\end{corollary}
That is, all characters giving isomorphic structures on a group algebra differ by the action of the automorphism $f$ on the groupoid $\Gamma$.

\begin{example}
In a group algebra $\mathbb {C}[G]\cong \bigoplus A^{i}$ such that $A^{1} \neq 0$, and given $a \in Z(G) \cap A^{1}$, the inner derivation generated by the element $a$ according to the formula from Proposition~\ref{Inn} defines a DG-algebra structure on $\mathbb {C}[G]$.
Note that the identity element of the group lies in $A^{0}$.
\begin{proof}
We have to show that $d^{2} =0$ for such a derivation. We know from Theorem \ref{DG_iff} that this condition is equivalent to 
\begin{center}
$\sum\limits_{h \in G}\chi(h, g)\chi(x,h) = 0$; $\forall g, x \in G.$
\end{center}
Consider the product $\chi(h, g)\chi(x,h)$. By the definition of an inner derivation, this product is not equal to 0 if and only if $(h,g)$ and $(x,h)$ $\in Hom(a,-a)$ or $Hom(-a,a)$.

Let $s(h,g) = s(x,h) = a$ and $t(h,g) = t(x,h) = -a$, then $(-1)^{|g|}g^{-1}h = (-1)^{|h|}h^{-1}x$ and $hg^{-1} = xh^{-1}$; but this is impossible because $a \in A^{1} \Rightarrow |h| \neq |g|.$

Let $s(h,g) = t(x,h) = a$ and $t(h,g) = s(x,h) = -a$, then $(-1)^{|g|}g^{-1}h = xh^{-1}$ and $hg^{-1} = (-1)^{|h|}h^{-1}x$; but this is similarly impossible because $|h| \neq |g|.$

The remaining two cases are similar.
\end{proof}
\end{example}

\begin{proposition}
\label{cen_der_th}
Let $d^{\tau}_{z}(g) = \tau(g)gz$ be the central derivation given by the element $z \in A^{1}$. It defines a DG-algebra structure if and only if $\tau(g)\tau(gz) = 0 \;\forall g \in G$.
\end{proposition}
\begin{proof}
The fact that $d^{\tau}_{z}(g) : A^{i} \rightarrow A^{i+1}$ follows from the fact that $z \in A^{1}$.\\
It remains to check under which conditions $d^{\tau}_{z}(g)^{2} = 0$.  To do this, we simply write by definition:
\begin{center}
$(d^{\tau}_{z})^{2}(g) = d^{\tau}_{z}(\tau(g)gz) = \tau(g)\tau(gz)gz^{2} \Rightarrow (d^{\tau}_{z})^{2}(g) = 0 \Leftrightarrow \tau(g)\tau(gz) = 0.$
\end{center}
\end{proof}

The example below illustrates a non-trivial condition that is implied by the previous Proposition.

\begin{example}\label{example3.2}
Let $N \triangleleft G$ be a normal subgroup in $G$ and $\exists f:\faktor{G}{N} \rightarrow \mathbb{Z}$ an epimorphism. Then $\mathbb{C}[G] \cong \bigoplus_{k \in \mathbb{Z}}  \mathbb{C}[\left\langle{w^{k}_{z}}\right\rangle N]$, where $z \in Z\left(\faktor{G}{N}\right)$ is such that $f(z) = 1$, and $\left\langle{w^{k}_{g}}\right\rangle$ is the set of words in which the total degree of  $z$ is $k$. In addition, the central derivation $d^{\tau}_{z}$, where 
$$
\tau(z^{k}g) = \begin{cases}
1 , \quad  k = 2n+1,\\
0 ,  \quad  k =2n,\\
\end{cases}
$$
defines a  DG-algebra structure on $\mathbb{C}[G]$.
\end{example}
Here $Z\left(\faktor{G}{N}\right)$ denotes the center of the group $\faktor{G}{N}$.

\begin{proof}
It is only necessary to check that the graded Leibniz rule is satisfied. To do this, we can make sure that the given $\tau$ is indeed a group character and then take advantage of Proposition \ref{cen_der_th}. It is not difficult to do this by simply going through all the possibilities.
\end{proof}
\begin{corollary}
In the graded group algebra $\mathbb {C}[\mathbb{Z}]\cong \bigoplus \mathbb {C}[x^k]$, where $\mathbb{Z} = \left\langle {x} \right\rangle$, it is possible to define a DG-algebra structure using the central derivation $d^{\tau}_x$, where 
$$
\tau(x^k) = \begin{cases}
1 , \quad  k = 2n+1,\\
0 ,  \quad  k =2n.\\
\end{cases}
$$
\end{corollary}

\begin{corollary}
Consider the Heisenberg integer group of 3×3 upper-triangular matrices
\begin{center}
 $H_{3} = \left\langle {x,y,z \mid z = xyx^{-1}y^{-1}, xz = zx, yz = zy} \right\rangle$.
\end{center}
It is well known that any element of this group can be represented as $y^{b}z^{c}x^{a}$. Then $f(y^{b}z^{c}x^{a}) = c \in \mathbb {Z}$ is an epimorphism into a group of integers. From here we get the grading $\mathbb {C}[H_{3}]\cong \bigoplus_{k \in \mathbb{Z}} \mathbb {C}[z^k\left\langle {x,y} \right\rangle]$. And the central derivation $d^{\tau}_{z}$, constructed similarly to Example \ref{example3.2}, gives a DG-algebra structure on $\mathbb{C}[H_{3}]$.
\end{corollary}

\subsubsection*{Acknowledgemets}
The results of Section \ref{sec-grouppoid} were obtained by A.~Arutyunov, the other results were obtained jointly. The work of A.~Arutyunov was supported by the RSF grant No. 25-11-00018.

\EditInfo{May 30, 2025}{September 19, 2025}{Ivan Kaygorodov  and David Towers}


{\small

\begin{thebibliography}{10}

\bibitem{Alekseev_2023}
A.~Alekseev, A.~Arutyunov, and S.~Silvestrov.
\newblock On $(\sigma ,\tau )$-derivations of group algebra as category
  characters.
\newblock {\em Springer Proceedings in Mathematics and Statistics}, 426:81--99,
  2023.

\bibitem{Arutyunov_2020}
A.~A. Arutyunov.
\newblock Derivation algebra in noncommutative group algebras.
\newblock {\em Proceedings of the Steklov Institute of Mathematics},
  308(1):22--34, 2020.

\bibitem{Arutyunov_2019}
A.~A. Arutyunov and A.~S. Mishchenko.
\newblock A smooth version of johnson's problem on derivations of group
  algebras.
\newblock {\em Sbornik: Mathematics}, 210(6):756--782, 2019.

\bibitem{arutyunov2017derivationsgroupalgebras}
A.~A. Arutyunov, A.~S. Mishchenko, and A.~I. Shtern.
\newblock Derivations of group algebras.
\newblock {\em Journal of Mathematical Sciences}, 248:709--718, 2020.

\bibitem{F_lix_2001}
Y.~F{\'e}lix, S.~Halperin, and J.-C. Thomas.
\newblock {\em Rational Homotopy Theory}.
\newblock Springer New York, 2001.

\bibitem{Lyndon_2001}
R.~C. Lyndon and P.~E. Schupp.
\newblock {\em Combinatorial Group Theory}.
\newblock Springer Berlin, Heidelberg, 2001.

\bibitem{Mao_2011}
X.~Mao.
\newblock Dg algebra structures on as-regular algebras of dimension 2.
\newblock {\em Science China Mathematics}, 54(10):2235--2248, 2011.

\bibitem{MAO2022389}
X.-F. Mao, X.-T. Wang, and M.-Y. Zhang.
\newblock \text{DG} algebra structures on the quantum affine $n$-space
  $\mathcal{O}_{-1}(k^n)$.
\newblock {\em Journal of Algebra}, 594:389--482, 2022.

\bibitem{Meirenken_2003}
E.~Meirenken.
\newblock Group actions on manifolds.
\newblock {\em Lecture Notes, University of Toronto},
  (math.toronto.edu/mein/teaching/LectureNotes/action.pdf), 2003.

\bibitem{Orlov_2023}
D.~Orlov.
\newblock Smooth \text{DG} algebras and twisted tensor product.
\newblock {\em Russian Mathematical Surveys}, 78(5):853–880, 2023.

\end{thebibliography}}
\end{document}